\documentclass[10pt,british]{article}

\usepackage[T1]{fontenc}
\usepackage[latin9]{inputenc}
\usepackage[a4paper]{geometry}
\usepackage{color}
\usepackage{babel}
\usepackage{mathrsfs}
\usepackage{amsmath}
\usepackage{amsthm}
\usepackage{amssymb}
\usepackage{setspace}
\usepackage[authoryear]{natbib}
\usepackage[unicode=true,pdfusetitle,
 bookmarks=true,bookmarksnumbered=false,bookmarksopen=false,
 breaklinks=false,pdfborder={0 0 1},backref=false,colorlinks=false]
 {hyperref}

\makeatletter
\theoremstyle{remark}

\theoremstyle{plain}

\theoremstyle{plain}

\theoremstyle{plain}

\theoremstyle{definition}

\usepackage{lscape}
\usepackage{changepage}
\date{}

\makeatother

\providecommand{\assumptionname}{Assumption}
\providecommand{\examplename}{Example}
\providecommand{\lemmaname}{Lemma}
\providecommand{\remarkname}{Remark}
\providecommand{\theoremname}{Theorem}

\begin{document}
\title{All-Pass Functions for Mirroring Pairs of Complex-Conjugated Roots of Rational Matrix Functions}
\author{Wolfgang Scherrer$^{a}$ and Bernd Funovits$^{a,b}$}

\maketitle
\thispagestyle{empty}

\section*{Proposed Running Head}

Mirroring Complex-Conjugated Roots

\section*{Affiliation}

\begin{singlespace}
	\textbf{$^{a}$TU Wien}
	
	Institute of Statistics and Mathematical Methods in Economics
	
	Econometrics and System Theory
	
	Wiedner Hauptstr. 8
	
	A-1040 Vienna
\end{singlespace}

and

\begin{singlespace}
	\textbf{$^{b}$University of Helsinki}
	
	Faculty of Social Sciences
	
	Discipline of Economics
	
	P. O. Box 17 (Arkadiankatu7)
	
	FIN-00014 University of Helsinki
\end{singlespace}

\section*{E-mail}

bernd.funovits@helsinki.fi

\pagebreak{}

\thispagestyle{empty}

\section*{Abstract}

We construct rational all-pass matrix functions with real-valued coefficients for mirroring pairs of complex-conjugated determinantal roots of a rational matrix.
This problem appears, for example, when proving the spectral factorization theorem, or, more recently, in the literature on possibly non-invertible or possibly non-causal vector autoregressive moving average (VARMA) models.
In general, it is not obvious whether the all-pass matrix function (and as a consequence the all-pass transformed rational matrix with initally real-valued coefficients) which mirrors complex-conjugated roots at the unit circle has real-valued coefficients.
Naive constructions result in all-pass functions with complex-valued coefficients which implies that the real-valued parameter space (usually relevant for estimation) is left.

Keywords: All-pass functions, Blaschke matrices, spectral factorization

JEL classification: C32, C50

\pagebreak{}

\section{Introduction}

\setcounter{page}{1}

It is well known that there are multiple spectral factors which generate the same spectral density. %
In the classical time series literature, \cite{Rozanov67} and \cite{Hannan70} use all-pass functions (also known as Blaschke factors) to mirror determinantal roots of spectral factors from inside to outside the unit circle when proving the spectral factorization theorem (for rational spectral densities).
More recently, all-pass functions play an important role in the literature on possibly non-invertible or possibly non-causal vector autoregressive moving average (VARMA) models, see \cite{LanneSaikkonen13}, \cite{funovits2020identifiability}.

It is not obvious whether the all-pass function (and as a consequence the all-pass transformed polynomial or rational matrix with initially real-valued coefficients) which mirrors complex-conjugated roots at the unit circle is real-valued and, to the best of our knowledge, there is no proof available in the literature which addresses this issue.
Naive constructions, as in \cite{GourierouxMR_svarma19}, result in all-pass functions with complex-valued coefficients which implies that the real-valued parameter space (usually relevant for estimation) is left.

Here, we show how to obtain all-pass functions with real-valued coefficients for mirroring pairs of complex-conjugated determinantal roots of a rational matrix function at the unit circle in three ways.
All approaches have the fact that they use a (static) orthogonal transformation to reduce the dimensionality of the problem in common.
More precisely, they start from the QR decomposition of the real and imaginary part of a normalized vector in the (right-) kernel of the rational matrix $k(z)$ which is evaluated at a determinantal zero with non-trivial imaginary part and whose complex-conjugated zeros are to be mirrored.
One approach parametrises consecutively unitary matrices and Blaschke factors in terms of the matrix $R$ of the QR decomposition and the real and imaginary part of the corresponding determinantal root of $k(z)$. 
While this approach is similar to the case of mirroring real-valued roots, it leaves the real-valued parameter space in intermediary steps and only at the end is it ensured that the coefficient matrices are indeed real-valued.
The other two approaches do not leave the real parameter space in any of the intermediate steps.
The second and more elegant approach uses the defining characteristics of all-pass functions more prominently and the calculations are less tedious.
The last construction is based on state space methods.

The remainder of this article is structured as follows: 
In Section 2, we define all-pass functions, including some special instances that will appear in our derivations, and review the spectral factorization problem.
Section 3 is of preparatory nature. 
We discuss how to reduce the problem of (cross-sectional) dimension $ n$ to either a $ 2$- or a $ 1$-dimensional problem (in the case of a pair of complex-conjugated roots or one real-valued root, respectively). %
Moreover, we deal with the case of a real-valued root and a degenerate case in order to be able to focus on the essential problem in the sections which follow.
Sections 4 to 6 deal with different approaches for mirroring complex-conjugated determinantal roots of a $ \left( 2 \times 2 \right) $ polynomial matrix.
In Section 4, we parametrise unitary matrices and Blaschke matrices in terms of the elements of $R$ and the real and imaginary part of the determinantal root such that their product has real coefficients and mirrors the given pair of complex-conjugated roots at the unit circle.
This procedure is closest to what can be found in existing literature.
In Section 5, we present a more elegant approach for mirroring complex-conjugated roots of a $ \left( 2 \times 2 \right) $ polynomial matrix. 
In contrast to the approach in Section 4, the real parameter space is not left in intermediate steps.
In Section 6, we discuss a state space approach for constructing all-pass functions with real-valued coefficients.

We use $z$ as a complex variable %
and define $i=\sqrt{-1}$. 
The transpose of an $\left(m\times n\right)$-dimensional matrix $A$ is represented by $A'$.
For the sub-matrix of $A$ consisting of rows $m_{1}$ to $m_{2}$, $0\leq m_{1}\leq m_{2}\leq m$, we write $A_{\left[m_{1}:m_{2},\bullet\right]}$ and analogously $A_{\left[\bullet,n_{1}:n_{2}\right]}$ for the sub-matrix of $A$ consisting of columns $n_{1}$ to $n_{2}$, $0\leq n_{1}\leq n_{2}\leq n$. 
The $n$-dimensional identity matrix is denoted by $I_{n}$, an $n$-dimensional diagonal matrix with diagonal elements $\left(a_{1},\ldots,a_{n}\right)$ is denoted by $\text{diag}\left(a_{1},\ldots,a_{n}\right)$. %

\section{All-Pass Functions, Blaschke Matrices, and the Rational Spectral Factorization Problem}

A multivariate rational all-pass function is an $(n\times n)$-dimensional matrix $V(z)$ whose entries are rational functions and which satisfies $V(z)V^{*}\left(\frac{1}{z}\right)=V^{*}\left(\frac{1}{z}\right)V(z)=I_{n}$.
The superscript asterisk is defined as $ m^{*}(z) = \overline{m\left(\bar{z}\right)}' $ such that the coefficient matrices of, say, a polynomial matrix are transposed and conjugated but the variable $ z $ remains unaffected.
The classical spectral factorization problem for rational spectral densities of full rank consists in factorizing an $ \left( n \times n \right) $-dimensional matrix $ f(z) $ whose elements are rational functions in $ z $ with real coefficients which satisfies $f(z) >0$ for all $|z|=1$ and $ f(z) = f^*\left(\frac{1}{z}\right)$ as $ f(z) = k(z) k^*\left(\frac{1}{z}\right) $. 
This factorization is not unique but a unique canonical factor $k(z)$, which is a rational function without poles inside or on the unit circle and without zeros inside the unit circle and $ k(0) $ is a lower-triangular matrix with positive diagonal elements, does exist.
This theorem is proved, for example, in \citet[page 47]{Rozanov67} and \citet[page 66]{Hannan70}.
While the canonical spectral factor $ k(z) $ described above is unique, it is easy to see that post-multiplying an all-pass function on $ k(z) $ results in a different spectral factor 
$ w(z) = k(z) V(z) $ which also satisfies $ f(z) = w(z) w^*\left(\frac{1}{z}\right)$. 
Furthermore, the all-pass functions $ V(z) $ can be chosen such that they do not introduce additional zeros or poles but rather mirror zeros at the unit circle, i.e. replace zeros $z$ by $1/z$.

As alluded to in the introduction, it is non-trivial to ensure that $ w(z) $ has real-valued coefficients even if pairs of complex-conjugated roots are mirrored together.
To fix ideas, consider a square polynomial matrix $ p(z) $ of dimension $ n $ and degree $ q $ with real-valued coefficients.
If $ p(z) $ has a total of $ m_z \leq n \cdot q$ determinantal roots, there are $ 2^{m_r + m_c} -1$ non-trivial all-pass functions such that $ w(z) $ again has real-valued coefficients, where $ m_r $ is the number of real-valued roots and $ 2 \cdot m_c $ is the number of complex-valued roots such that $ m_r + 2\cdot m_c = m_z $.

Before we start with the main part, we introduce some special all-pass functions which will be useful in constructing real-valued Blaschke matrices which mirror complex-conjugated roots.
An elementary Blaschke factor at $\alpha \in \mathbb{C}$ is of the form\footnote{
Sometimes, the Blaschke factor is defined with an additional factor $\alpha/ |\alpha|$.
However, this factor is not well defined if $\alpha = 0$.} $B(z, \alpha) = \frac{1-\bar{\alpha} z}{-\alpha + z}$.
A squared Blaschke factor pertaining to the pair of complex-conjugated roots $\alpha_{\pm}=\alpha_r \pm i\alpha_i$, where $ \alpha_i > 0 $ and $ \alpha_- = \overline{\alpha_{+}} $, is defined as $B_{sq}(z, \alpha_{\pm}) = \frac{1-\alpha_{-} z}{-\alpha_{+} + z} \frac{1-\alpha_{+} z}{-\alpha_{-} + z} $ and thus has real-valued coefficients.
Lastly, a bivariate Blaschke factor pertaining to $\alpha_{\pm}$ and the non-zero vector $w \in \mathbb{C}^{2\times 1}$ is given as $B_2(z, \alpha_{\pm}, w) = a^{-1}(z) b(z)$, where $a(z) = \left( z-\alpha_+ \right) \left( z-\alpha_- \right) $ (which has real-valued coefficients) and $b(z)$ is a $(2 \times 2)$ polynomial matrix with highest degree $2$ and which is of reduced rank at $z=\alpha_{+}^{-1}$, $z=\alpha_{-}^{-1}$, $z=\alpha_{+}$ and $z=\alpha_{-}$.
We will construct $b(z)$ such that its coefficients are real-valued,  the column space of $b(\alpha_+)$ is spanned by a given (non-trivial) vector $w \in \mathbb{C}^{2 \times 1}$, and $B_2(z, \alpha_{\pm}, w)$ is all-pass.
The resulting all-pass function is, of course, not unique since it could, for example, be post-multiplied by an orthogonal matrix.

\section{Reducing the Problem to the Essential}\label{sec:reduce}

Any rational matrix $k(z)$ may be factored as $k(z)=q^{-1}(z) p(z)$ with two polynomial matrices $q(z), p(z)$ such that the zeros of $k(z)$ are the the zeros of $p(z)$. 
Hence, it suffices to consider polynomial matrices for which we treat three different cases:
The case of a real-valued determinantal root of $ p(z) $, the degenerate case of a pair of complex-conjugated roots which is almost identical to the real-valued case, and finally the non-degenerate case (which will be the subject of subsequent sections).
We will use transformations of the form 
\begin{equation}\label{eq:blaschke_reduce_a}
	\tilde{p}(z) = p(z) Q \begin{pmatrix}
		V\left( z \right) & 0 \\
		0 & I_m
	\end{pmatrix},
\end{equation}
where $Q$ is a real-valued orthogonal matrix, $ V(z) $ is an all-pass function (with real-valued coefficients) corresponding to one of $ B(z, \alpha) $, $B_{sq}(z, \alpha_{\pm})$, or $B_2(z, \alpha_{\pm}, w)$ defined in the previous section, and $ I_m $ an identity matrix of appropriate dimension. %
By construction, $\tilde{p}(z) \tilde{p}^{*}\left(\frac{1}{z}\right) = p(z) p^{*}\left(\frac{1}{z}\right)$ holds. 
However, we have to make sure  that $\tilde{p}(z)$ is a polynomial matrix with real-valued coefficients.

First, let us consider the case of a real-valued zero $\alpha_r \in \mathbb{R}$ of $ \det\left( p\left( z \right) \right) $ for which we may choose $v\in\mathbb{R}^{n \times 1}$ in the right-kernel of $ \det\left( p\left( \alpha_r \right) \right) $. 
All entries of $p(z)v$ have $\left( z-\alpha_r \right)$ as a common factor which implies that $p(z)vB\left( z,\alpha_r \right)$ is a polynomial vector which is equal to zero when evaluated at $z=\alpha_r^{-1}$. 
Thus, we set $V(z)=B\left( z,\alpha_r \right)$ in \eqref{eq:blaschke_reduce_a} where $Q$ is an orthogonal matrix whose first column is equal to the normalized vector $v/\|v\|$.

Secondly, we consider a pair of complex-conjugated roots $ \alpha_{\pm} $ for which a non-zero vector\footnote{
	This case is equivalent to the condition that there exists a vector $v\in\mathbb{C}^n$, $p\left( \alpha_+ \right)v=0$ such that $v$ and $\overline{v}$ are linearly dependent in $\mathbb{C}^n$. 
	This situation occurs for $n=1$ for all complex-conjugated zeros. 
	However, this case is, in a certain sense, non-generic for $ n > 1 $. 
} 
$v\in\mathbb{R}^n$ such that $p(\alpha_+)v=0$ exists. 
Since $p$ has real-valued coefficients, the complex-conjugate $\alpha_- = \overline{\alpha_+}$ is also a determinantal zero.
Therefore, $p\left( \alpha_+ \right)v=p\left( \alpha_- \right)v=0$ and hence the elements of $p(z)v$ have $\left( z - \alpha_+ \right) \left( z - \alpha_- \right)$ as a common factor. 
Thus, we set $V(z)=B_{sq}\left( z,\alpha_{\pm} \right)$ in \eqref{eq:blaschke_reduce_a} where $Q$ again is an orthogonal matrix with the first column equal to the normalized vector $v/\|v\|$. 

Last, we treat the case where the right-kernel of $p\left( \alpha_+ \right)$ does not contain a real non-zero vector. 
Thus, for any vector $v\in\mathbb{C}^{n \times 1}$, $\left( \alpha_+ \right)v=0$ the two vectors $v,\overline{v}$ are linearly independent. 
There exists a real-valued semi-orthogonal matrix\footnote{
	Consider the QR decomposition $(v_r, v_i)=Q_1 R$, where $v=v_r + iv_i$, $v_r, v_i \in \mathbb{R}^{n \times 1}$, and set $w=Q_1'v = R \left(\begin{smallmatrix}1 \\ i\end{smallmatrix}\right)$.} 
$Q_1\in\mathbb{R}^{n\times 2}$, $Q_1'Q_1 = I_2$ and a vector $w\in\mathbb{C}^{2 \times 1}$  such that $v=Q_1w$ and $\overline{v}=Q_1 \overline{w}$. 
We expand $Q_1$ to an orthogonal matrix $Q=(Q_1,Q_2)$ and use $V(z)=B_2\left( z,\alpha_\pm,w \right) = b(z) \left[ \left( z-\alpha_+ \right) \left( z-\alpha_{-} \right) \right]^{-1}$ in the transformation~\eqref{eq:blaschke_reduce_a}.
Note that the polynomial matrix $b(z)$ is constructed such that $b(\alpha_+)=wm'$, for a vector $m\in\mathbb{C}^{2 \times 1} $. 
Therefore, $p\left( \alpha_+ \right) Q_1 b\left( \alpha_+ \right) =p\left( \alpha_+ \right) Q_1 w m' =p\left( \alpha_+ \right) v m' = 0\in \mathbb{C}^{n\times 2}$ which shows that all entries of $p(z) Q_1 b(z)$ have $\left(z- \alpha_+ \right)$ as a common factor. 
Analogously, it follows that $\left( z-\alpha_- \right)$ is a common factor of the elements in $p(z)Q_1 b(z)$. 
Together, this implies that $p(z)Q_1 B_2(z,\alpha_\pm,w)$ is in fact polynomial.
In the following three sections, we present three alternative constructions for the bivariate Blaschke matrix $B_2(z,\alpha_\pm,w)$.

\section{Parametrising All-Pass Functions Consecutively}

Here, we construct the all-pass function $ V(z) $ in the  non-degenerate case of complex-conjugated determinantal roots $ \alpha_{\pm} $ explicitly.
We parametrise unitary $(2\times2)$-dimensional matrices $V_{\beta}$, $V_{\gamma}$, and $V_{\delta}$ in terms of the parameters in $ R $ (obtained from the QR decomposition $ \left( v_r, v_i \right) = Q_1 R $ from the previous section) and $ \alpha_{\pm} $ such that 
\begin{equation}\label{eq:blaschke_successive}
	 V_{\beta} 
	\left( \begin{smallmatrix}
		B(z, \alpha_+) & 0 \\
		0 &  1
	\end{smallmatrix} \right)  V_{\gamma} 
	\left( \begin{smallmatrix}
		B(z, \alpha_-) & 0 \\
		0 &  1
	\end{smallmatrix} \right)
	V_{\delta}
\end{equation}
has real-valued coefficient matrices.
We parametrise the  $(2\times 2)$ unitary matrices in this section with two parameters $\left( \phi_1, \phi_2 \right)$ as 
$\left( \begin{smallmatrix}  \cos(\phi_1) e^{i \phi_2} & - \sin(\phi_1) \\\sin(\phi_1) & \cos(\phi_1) e^{-i \phi_2}  \end{smallmatrix} \right)$.

First, we obtain $V_{\beta}$ by choosing $\beta_1, \beta_2$ such 
that $R\left(\begin{smallmatrix} 1 \\ i\end{smallmatrix}\right)$ is in the span of
$\left(\begin{smallmatrix} \cos(\beta_1) e^{i \beta_2} \\  \sin(\beta_1) \end{smallmatrix}\right)$.
More specifically, for $R = \left( \begin{smallmatrix} a & b \\ 0 & c \end{smallmatrix} \right)$, where $a$ and $c$ are positive and $a^2+b^2+c^2=1$ by construction (as a consequence of the QR decomposition and the fact that $\| v_r \|^2 + \| v_i \|^2 = 1$), 
the column space spanned by\footnote{
	Remember that for $z=x+iy$, the polar representation $z=r\cdot\cos(\phi) + i\cdot r\cdot\sin(\phi)$ can be obtained with $r = \sqrt{x^2+y^2}$ and, for $r>0$, $\phi = \arccos\left(\frac{x}{r}\right)$ when $y >0$ and $\phi = -\arccos\left(\frac{x}{r}\right)$ when $y\leq0$.
	Note that $\frac{a}{c}$ is always positive by construction.
} 
$R\left(\begin{smallmatrix} 1 \\ i\end{smallmatrix}\right)$ is equal to the one spanned by $\left(\begin{smallmatrix} \frac{b}{c} - i \frac{a}{c} \\ 1 \end{smallmatrix}\right) = 
\left(\begin{smallmatrix} 
	\frac{\sqrt{a^2+b^2}}{c} \cdot \exp\left(-i\cdot \arccos\left(\frac{b}{\sqrt{a^2+b^2}}\right) \right) \\ 
	1 
\end{smallmatrix}\right)$.
Normalising this vector leads to 
$$
V_{\beta} = 
\left( \begin{smallmatrix}
	\sqrt{a^2+b^2} \cdot e^{-i\cdot \arccos\left(\frac{b}{\sqrt{a^2+b^2}} \right)} & - c \\
	c & 	\sqrt{a^2+b^2} \cdot e^{i\cdot \arccos\left(\frac{b}{\sqrt{a^2+b^2}} \right)}
\end{smallmatrix} \right)
=
\left( \begin{smallmatrix} 
	\cos(\beta_1) e^{i \beta_2} & -\sin(\beta_1)  \\  
	\sin(\beta_1) & \cos(\beta_1) e^{-i \beta_2}
\end{smallmatrix} \right).
$$

Secondly, $V_{\gamma}$ is determined by setting $\gamma_1, \gamma_2$ (which depend in turn on $\beta_1, \beta_2$, and $\alpha_{\pm}$) such that $R\left(\begin{smallmatrix}1 \\ -i\end{smallmatrix}\right)$, which is in the span of $\overline{V_{\beta, [\bullet, 1]}} = \left(\begin{smallmatrix} \cos(\beta_1) e^{-i \beta_2} \\ \sin(\beta_1) \end{smallmatrix}\right)$, is in turn in the span of 
$V_{\beta}  \left( \begin{smallmatrix} B(\alpha_-,\alpha_+) & 0 \\ 0 & 1 \end{smallmatrix} \right)  V_{\gamma,[\bullet,1]}$. 
Thus, the first column of $V_{\gamma}$ should be equal to a normalised version of $\left( \begin{smallmatrix} B(\alpha_-,\alpha_+)^{-1} & 0 \\ 0 & 1 \end{smallmatrix} \right)  V_{\beta}^{-1}\left(\begin{smallmatrix} \cos(\beta_1) e^{-i \beta_2} \\ \sin(\beta_1) \end{smallmatrix}\right)$ which eventually leads to the first column of $V_{\gamma}$ being equal to
$
\frac{1}{k}
\left( \begin{smallmatrix} 
	B(\alpha_-,\alpha_+)^{-1} \left( \cos\left(\beta_2\right)^2 \cdot e^{-i2\beta_1} + \sin\left(\beta_2\right)^2  \right) \\  
	2i\sin(\beta_2) \cos(\beta_2) \sin(\beta_1)
\end{smallmatrix} \right),
$
where $k$ is a normalising constant.

Last, $V_{\delta}$ is the inverse of $V_{\beta} \left( \begin{smallmatrix}  B(1,\alpha_+) & 0 \\ 0 & 1 \end{smallmatrix} \right)  V_{\gamma}  \left( \begin{smallmatrix}  B(1,\alpha_-) & 0 \\ 0 & 1 \end{smallmatrix} \right)$ and again a function of $\beta_1, \beta_2$, and $\alpha_{\pm}$.
Straight-forward computation verifies that the coefficient matrices in $$V_{\beta}  
\left(
\begin{smallmatrix}
	B(z, \alpha_+) & 0 \\
	0 &  1
\end{smallmatrix}
\right)  \left( -\alpha_+ + z \right) \left( -\alpha_- + z \right)  V_{\gamma}  
\left(
\begin{smallmatrix}
	B(z, \alpha_-) & 0 \\
	0 &  1
\end{smallmatrix}
\right)  \left( -\alpha_+ + z \right) \left( -\alpha_- + z \right)  V_{\delta},
$$ which is equal to
$$
V_{\beta}  \left( \begin{smallmatrix} \left(-\alpha_- +z\right) \left(1-\alpha_- z\right) & 0 \\ 0 & \left(-\alpha_+ +z\right) \left(-\alpha_- +z\right) \end{smallmatrix} \right)  
V_{\gamma}  \left( \begin{smallmatrix} \left(-\alpha_+ +z\right) \left(1-\alpha_+ z\right) & 0 \\ 0 & \left(-\alpha_+ +z\right) \left(-\alpha_- +z\right) \end{smallmatrix} \right)  
V_{\delta},
$$
are indeed real.
Since $  \left( z-\alpha_+ \right) \left( z-\alpha_-\right)$ has real-valued coefficients, it follows that the polynomial matrix in \eqref{eq:blaschke_successive} has real-valued coefficients as well.

\section{Polynomial Approach}

While the approach in the previous section is straight-forward and similar to the real-valued case in the sense that one root is mirrored after the other, it is tedious and does not create additional insights into the structure of the problem.
Here, we present a more elegant construction which uses the defining characteristics of all-pass functions in a more prominent way.

We proceed in the following steps.
Firstly, for a given non-singular matrix $A \in \mathbb{R}^{2 \times 2} $, we construct an all-pass function 
of the form $V(z) = \left( I_2 - A z \right)^{-1} \left( I_2 - B z \right)  T^{-1}$. We assume that $A$ does 
not have eigenvalues on the complex unit circle, i.e. the moduli of the eigenvalues are \emph{not} equal to one. 
In the second step, we show how to choose the matrix $A$ such that $V(z)=B_2(z, \alpha_\pm, w)$. 

We start with the case where the eigenvalues of $A$ are inside the unit circle. 
The matrix $ V(z) = (I-Az)^{-1}(I-Bz)T^{-1} $ is all pass if and only if 
\begin{equation}
\label{eq:all_pass_eq}
	(T')^{-1}(I-B'z^{-1})(I-A'z^{-1})^{-1}(I-Az)^{-1}(I-Bz)T^{-1} = I_2
\end{equation}
holds.
In turn, this is equivalent to 
\begin{align*}
(I-A'z^{-1})^{-1} (I-Az)^{-1} &= (I-B'z^{-1})^{-1}T'T(I-Bz)^{-1} \\
\left(\sum_{k=0}^\infty (A')^k z^{-k} \right)\left(\sum_{k=0}^\infty A^k z^k \right) &= 
\left(\sum_{k=0}^{\infty} (B')^{-k} z^{k} \right) (B')^{-1}T'TB^{-1}\left(\sum_{k=0}^{\infty} B^{-k} z^{-k} \right) \\
\Gamma_0 + \sum_{k=1}^\infty (\Gamma_k z^k + \Gamma_k' z^{-k}) &= 
\Omega_0 + \sum_{k=0}^\infty (\Omega_k z^k + \Omega_k' z^{-k})
\end{align*}
with
\begin{align*}
\Gamma_0 &= \sum_{k=0}^\infty (A')^k A^k = A'\Gamma_0 A + I_2 \\
\Omega_0 &= \sum_{k=0}^\infty (B')^{-k} (B')^{-1}T'TB^{-1} B^{-k} = (B')^{-1}\Omega_0 B^{-1} + (B')^{-1}T'TB^{-1}\\
\Gamma_k &= \Gamma_0 A^k = (\Gamma_0 A \Gamma_0^{-1})^k \Gamma_0 \\
\Omega_k &= (B')^{-k} \Omega_0 
\end{align*}
In the above derivations, we have silently assumed that the eigenvalues of $B$ have moduli larger than one in order to obtain the Laurent series expansion of $(I-Bz)^{-1} $.
Now, it is straight-forward to see that \eqref{eq:all_pass_eq} holds if and only if $$B = \Gamma_0^{-1} (A')^{-1} \Gamma_0 \text{ and } T'T = B'\Gamma_0 B - \Gamma_0.$$
In particular, we see that $B$ indeed has eigenvalues outside the unit circle. 
Note also that $T$, which may, for example, be computed via a Cholesky decomposition of $B'\Gamma_0 B - \Gamma_0$, is non-singular. 

The case where the eigenvalues of $A$ are outside the unit circle may be solved analogously. 
In this case, we have to set  $B$ and $T$ as 
\begin{gather*}
B = \Gamma_0^{-1} (A')^{-1} \Gamma_0 \text{ and } T'T = \Gamma_0 - B' \Gamma_0 B \text{ where } \\ 
\Gamma_0 = (A')^{-1} \Gamma_0 A^{-1}  + (A')^{-1}  A^{-1} .
\end{gather*}

Now that we have shown how to construct $ \left( B, T \right) $ for given $ A $ such that $ V(z) = (I-Az)^{-1}(I-Bz)T^{-1} $ is all-pass, suppose that we are given a non-zero vector $w\in\mathbb{C}^{2 \times 1}$ (where $w$ and $\overline{w}$ are linearly independent) and set 
\[
A=
(w, \overline{w}) 
\begin{pmatrix}
	\alpha_+^{-1} & 0 \\
	0 & \alpha_-^{-1}
\end{pmatrix}
(w, \overline{w})^{-1} \in\mathbb{R}^{2\times 2}.
\]
Then, it may be easily seen that 
\begin{align*}
V(z) &= 
(I_2-Az)^{-1}  (I-Bz)T^{-1} \\
&= 
(\underbrace{\left( z-\alpha_+ \right) \left( z-\alpha_- \right) }_{=a(z)})^{-1}
\underbrace{
(w, \overline{w}) 
\begin{pmatrix}
	z - \alpha_- & 0 \\
	0 & z- \alpha_+
\end{pmatrix}
(w, \overline{w})^{-1} 
(I-Bz)T^{-1}
}_{=b(z)} = a^{-1}(z)b(z).
\end{align*}
It is immediate that the column space of $b\left( \alpha_+ \right)$ is spanned by $w$, which proves $V(z)=B_2(z,\alpha_\pm, w)$.

\section{State Space Construction}

Another elegant approach to construct the bivariate Blaschke factor $B_2(z, \alpha_{\pm}, w)$ is based on state space methods. 
Any rational matrix function $k(z)$ (with real valued coefficients) which has no pole at $z=0$ may be represented as 
$$
k(z) = C(z^{-1}I -A)^{-1}B + D
$$
with suitably chosen matrices $(A,B,C,D)$. 
This quadruple of matrices is called a state space realization of the rational function $k(z)$. 
Here, we give a construction for a state space realization for $B_2(z, \alpha_{\pm}, w) = a^{-1}(z)b(z) = C(z^{-1}I -A)^{-1}B + D$.

\subsection{Fixing the Poles of the All-Pass Function: Determining A}

For a given pair of complex-conjugated zeros $\alpha_{\pm}$ we set
$A = \left( \begin{smallmatrix} \lambda_r & \lambda_i \\ -\lambda_i & \lambda_r \end{smallmatrix} \right)$, where 
$\lambda_+ = (\lambda_r + i\lambda_i) = \alpha_+^{-1}$. 
The rational matrix $B_2(z, \alpha_{\pm}, w)$ may then be factored as 
$$
\begin{aligned}
	B_2(z, \alpha_{\pm}, w) &= C \begin{pmatrix}
		z^{-1} - \lambda_r & -\lambda_i \\
		\lambda_i & z^{-1} - \lambda_r 
	\end{pmatrix}^{-1} B + D \\
	&= (z^{-2} - 2\lambda_r z^{-1} + \lambda_r^2+\lambda_i^2)^{-1} 
	C \begin{pmatrix} z^{-1} - \lambda_r & \lambda_i \\ -\lambda_i & z^{-1} - \lambda_r \end{pmatrix} B + D \\
	&= \underbrace{(1 - 2\lambda_r z + |\lambda|^2 z^2)^{-1}}_{=:a(z)^{-1}} 
	C \begin{pmatrix}  z - \lambda_r z^2 & \lambda_i z^2 \\ -\lambda_i z^2 & z - \lambda_r z^2 \end{pmatrix} B + D \\
	&= a(z)^{-1}
	\underbrace{\left( C \begin{pmatrix}  z - \lambda_r z^2 & \lambda_i z^2 \\ -\lambda_i z^2 & z - \lambda_r z^2 \end{pmatrix} B + a(z)D \right)}_{=:b(z)}  = a^{-1}(z) b(z)                                          
\end{aligned}.
$$
Note that $a(z)=(1-\alpha_+^{-1}z)(1-\alpha_{-}^{-1}z)$ which shows that $B_2(z, \alpha_{\pm}, w)$ 
has poles at $z=\alpha_+$ and $z=\alpha_{-}$.

\subsection{Fixing the Column-Space at $\alpha_{\pm}$: Determining C}

Next, we determine $C$ such that the column-space of $b\left(\alpha_+\right)$ is spanned by a given column 
vector $w \in \mathbb{C}^{2 \times 1}$, in our case 
$w = R\left(\begin{smallmatrix}1 \\ i\end{smallmatrix}\right)$.
Note that
$$
\begin{aligned}
	b(\alpha_{+}) & = \alpha_{+}^2 C \begin{pmatrix} \alpha_{+}^{-1} - \lambda_r & \lambda_i \\ 
		-\lambda_i & \alpha_{+}^{-1} - \lambda_r \end{pmatrix} B + 
	a(\alpha_{+}) D \\
	& = \alpha_{+}^2 C 
	\begin{pmatrix} (\lambda_r + i \lambda_i) - \lambda_r & \lambda_i \\ 
		-\lambda_i & (\lambda_r + i \lambda_i) - \lambda_r \end{pmatrix} B 
	= \alpha_{+}^2 C 
	\begin{pmatrix} i \lambda_i & \lambda_i \\ 
		-\lambda_i & i \lambda_i \end{pmatrix} B. 
\end{aligned}             
$$

Therefore, we set $C \left(\begin{smallmatrix}  i \lambda_i \\ -\lambda_i \end{smallmatrix}\right)  = \lambda_i \|w\|^{-1} w$, i.e. $C = \|w\|^{-1} \left(\begin{smallmatrix} w_i & -w_r \end{smallmatrix}\right)$ where $w_r, w_i$ denote the real and imaginary parts of $w$.

\subsection{Ensuring All-Pass Property: Determining B and D} 

Finally, we construct $B,D$ (for given $A,C$) 
such that the rational matrix $B_2(z, \alpha_{\pm}, w)$ is indeed all-pass.
In the following we will represent a state space realization $(A,B,C,D)$ 
by the matrix  
$
\left(\begin{smallmatrix}
	A  &|& B   \\
	\hline
	C  &|& D
\end{smallmatrix}\right)
$ 
because many operations may be elegantly represented in this notation.  In particular, 
the product $B'_2\left(\frac{1}{z}, \alpha_{\pm}, w\right) B_2(z, \alpha_{\pm}, w)$ has a realization\footnote{
	In general, the multiplication of two rational functions $k_1(z)$ and $k_2(z)$ of appropriate 
	dimensions and parametrized as two state space systems $\left(\begin{array}{c|c}
		A_1  & B_1   \\
		\hline
		C_1  & D_1
	\end{array}\right)$ and 
	$\left(\begin{array}{c|c}
	  A_2 & B_2 \\ 
	  \hline 
	  C_2 & D_2
	  \end{array}\right)$ 
	such that 
	$k_1(z) \cdot k_2(z) = \left( C_1 (z^{-1}I_n - A_1)^{-1}B_1 + D_1 \right)\left( C_2 (z^{-1}I_n - A_2)^{-1}B_2 + D_2 \right)$ 
	results in the state space system 
	$$
	\left(\begin{array}{@{}cc|c@{}}
		A_1      & B_1 C_2 & B_1 D_2   \\
		0        & A_2     & B_2 \\  \hline
		C_1      & D_1 C_2 & D_1 D_2
	\end{array}\right).
	$$
	For a rational function $k(z)=C(z^{-1}I_n-A)^{-1}B+D$ where $A$ is non-singular,
	the function $k'\left(\frac{1}{z}\right)$ may be written in terms of $ z $ as 
	$$
	\begin{array}{rl}
		k'\left(\frac{1}{z}\right)
		&= B' \left(zI_n - A'\right)^{-1} C' + D'
		 = B' \left(\left(z(A')^{-1} - I_n\right)A'\right)^{-1} C' + D' \\
		&= \left[-B' (A')^{-1}\right] \left(I_n - z(A')^{-1}\right)^{-1} C' + D'
		= \left[-B' (A')^{-1}\right] \left( \sum_{j=0}^{\infty} (A')^{-j}z^j \right) C' + D' \\
		&= \left[-B' (A')^{-1}\right] \left( \sum_{j=1}^{\infty} A'^{-j}z^j \right) C' + \left( D' - B' (A')^{-1} C' \right) \\
		&= \left[-B' (A')^{-1}\right] \left(\frac{1}{z}-(A')^{-1}\right)^{-1} \left[(A')^{-1} C' \right] +
		  \left( D' - B' (A')^{-1} C' \right),
	\end{array}
	$$
	i.e. $k'(z^{-1})$ has a state space realization
	$$
	\left(\begin{array}{c|c}
		(A')^{-1} & (A')^{-1} C' \\
		\hline
		-B' (A')^{-1} & D' - B' (A')^{-1} C'
	\end{array}\right).
	$$
} %
given by
$$
\left(\begin{array}{@{}cc|c@{}}
	(A')^{-1}     & (A')^{-1} C' C                        & (A')^{-1} C' D          \\
	0           &  A                                  & B                     \\ \hline
	-B'(A')^{-1}  & \left(D'-B'(A')^{-1} C'\right) C      & D'D - B' (A')^{-1} C' D
\end{array}\right).
$$

A state transformation is a mapping involving a non-singular matrix $ M $ (of the same dimension as $A$) 
which maps one state space realization $ \left( A, B, C, D \right) $ to another one of the form 
$ \left( M A M^{-1}, M B, C M^{-1}, D \right) $.
Both generate the same transfer function.
Applying the state transformation $ \left( \begin{smallmatrix}
	I_n & X \\
	0 & I_n
\end{smallmatrix} \right) $ to the model above results in 
$$
\left(\begin{array}{@{}cc|c@{}}
		(A')^{-1}     & (A')^{-1} C' C + X A - (A')^{-1} X                      & (A')^{-1} C' D + X B    \\
		0           &  A                                                  & B                     \\ \hline
		-B'(A')^{-1}  & \left(D'-B'(A')^{-1} C'\right) C + B' (A')^{-1} X       & D'D - B' (A')^{-1} C' D
	\end{array}\right).
$$

We now choose $X$ and $B$ such that the (1,2), the (1,3) and the (3,2) blocks of this realization are zero.
That is,

$$
\begin{array}{rclcrcl}
	(A')^{-1} C' C + X A - (A')^{-1} X                &=& 0 & \Longleftrightarrow     & C' C + A' X A - X   &=& 0 \\
	(A')^{-1} C' D + X B                            &=& 0 & \Longleftrightarrow     & \left( X, (A')^{-1} C' \right) 
	     \begin{pmatrix} B \\ D \end{pmatrix}         &=& 0 \\
	\left(D'-B'(A')^{-1} C'\right) C + B' (A')^{-1} X &=& 0 & \Longleftrightarrow     & 
	     \left( X' A^{-1} - C' C A^{-1}, C' \right) \begin{pmatrix} B \\ D \end{pmatrix} &=& 0 \\
	& &   & \Longleftrightarrow     & \left( A' X', C' \right) \begin{pmatrix} B \\ D \end{pmatrix}  &=& 0 \\
\end{array}
$$

First, we obtain $X$ as solution of the Lyapunov equation obtained from block (1,2).
Second, we obtain $B$ as a function of $D$ from the equation obtained from block (1,3).
In particular, we set $B = -X^{-1} (A')^{-1} C' D$. 
By this choice, the (3,2) block is also made zero. 

The matrix $D$ is then constructed such that the (3,3) block is equal to the identity matrix, i.e. 
$
D'D - B' (A')^{-1} C' D = I_2.
$
Together with the above, we thus have
$$
\begin{array}{rl}
	I_2 &= D'D + D' C A^{-1} X^{-1} (A')^{-1} C' D \\
	&= D' \left( I_n +  C A^{-1} X^{-1} (A')^{-1} C' \right) D
\end{array}
$$
and we may obtain $D$ from a Cholesky decomposition of $\left( I_2 +  C A^{-1} X^{-1} (A')^{-1} C' \right)$.

It follows that for given $ \left( A, C \right) $ and subsequent determination of $ \left( B, D \right) $ 
involving a state transformation, we obtain that the product 
$$
B'_2\left(\frac{1}{z}, \alpha_{\pm}, w\right) B_2(z, \alpha_{\pm}, w) = 
\underbrace{
	\begin{pmatrix} -B'(A')^{-1} & 0 \end{pmatrix} 
	\begin{pmatrix}  I_n z^{-1} - (A')^{-1} & 0 \\ 0 & I_n z^{-1} - A \end{pmatrix}^{-1} 
	\begin{pmatrix} 0 \\ B \end{pmatrix}
}_{=0} 
+ I_2
$$
is indeed equal to the identity matrix for all $z\in\mathbb{C}$, i.e. $B_2(z,\alpha_\pm,w)$ is an all-pass function 
with real-valued coefficients.

This construction of all-pass rational matrix functions in state space form is especially useful if the rational matrix $ k(z) $ is in state space form and may be generalised in various ways. 
One may generalise the construction to $(n\times n)$ matrices with $n=1$ or $n>2$, it is possible to mirror several (real and/or pairs of complex-conjugated zeros) in one step, and one may also construct all-pass matrices which mirror poles analogously. 

\section{Acknowledgments}

Financial support by the Research Funds of the University of Helsinki as well as by funds of the Oesterreichische Nationalbank (Austrian Central Bank, Anniversary Fund, project number: 17646) is gratefully acknowledged. 

\section{Data Availability Statement}

There is no data involved in this study.

\end{document}